\documentclass{article}

\usepackage{amssymb}

\newcommand{\proof}{{\bf Proof:  }}
\newcommand{\remark}{{\bf Remark:  }}

\newcommand{\dimv}{\underline{\dim}}

\newcommand{\hb}{\newline\hspace*{\fill}$\Box$}

\newtheorem{theorem}{Theorem}[section]
\newtheorem{lemma}[theorem]{Lemma}
\newtheorem{definition}[theorem]{Definition}
\newtheorem{proposition}[theorem]{Proposition}
\newtheorem{corollary}[theorem]{Corollary}

\begin{document}
\parindent0pt
\title{\bf Poisson automorphisms and quiver moduli}

\author{Markus Reineke\\ Fachbereich C - Mathematik\\ Bergische Universit\"at Wuppertal\\ D - 42097 Wuppertal, Germany\\
e-mail: reineke@math.uni-wuppertal.de}
\date{}
\maketitle

\begin{abstract}
A factorization formula for certain automorphisms of a Poisson algebra associated with a quiver is proved, which involves framed versions of moduli spaces of quiver representations. This factorization formula is related to wall-crossing formulas for Donaldson-Thomas type invariants of M. Kontsevich and Y. Soibelman \cite{KS}.
\end{abstract}
\section{Introduction}

In \cite{KS}, a framework for the definition of Donaldson-Thomas type invariants for Calabi-Yau categories endowed with a stability structure is developed. One of the key features of this setup is a wall-crossing formula for these invariants, describing their behaviour under a change of stability structure in terms of a factorization formula for automorphisms of certain Poisson algebras defined using the Euler form of the category.\\[1ex]
In this paper, we study such factorization formulas using quiver representations, their moduli spaces, and Hall algebras (for different such approaches, see \cite{BTL,GPS}).\\[1ex]
With a quiver without oriented cycles, we associate a Poisson structure on a formal power series ring and study factorizations of a certain automorphism into infinite ordered products associated with stabilities for the quiver (see section \ref{statementresult} for precise definitions). Our main result, Theorem \ref{mainth}, describes these factorizations in terms of generating series for Euler characteristics of framed versions of moduli spaces of representations of the quiver, called smooth models in \cite{SM}. This result yields a weak form of an integrality conjecture for Donaldson-Thomas type invariants of \cite{KS}.\\[1ex]
We approach this theorem by producing the desired factorizations in the Hall algebra of the quiver, in terms of the Harder-Narasimhan recursion of \cite{HNS}; section \ref{hallhn} therefore adapts some of the material of \cite{HNS} to the present setup. Using an evaluation map, simple identities in the Hall algebra often lead to interesting identities in a skew formal power series ring associated with the quiver (see \cite{RICRA12} for some instances of this principle). In the present setup, the latter ring is viewed as a quantization of the Poisson algebra; we develop a simple algebraic setup for constructing automorphisms of the Poisson algebra from a class of formal series in section \ref{qd}. The key feature is that certain series in the skew formal power series ring induce non-trivial automorphisms of the Poisson algebra via conjugation. That this property applies to the generating series induced from the Harder-Narasimhan recursion follows from a formula for the Betti numbers of smooth models of \cite{SM}; this is explained in section \ref{smpt}. We end with two classes of examples in section \ref{examples}: considering generalized Kronecker quivers, we can prove a weak form of the integrality conjecture of \cite{KS}. We also relate the factorization formula Theorem \ref{mainth} to a factorization formula of \cite{GPS} in terms of Gromov-Witten theory. Finally, we consider Dynkin quivers, for which the relevant moduli spaces are trivial.\\[3ex]
{\bf Acknowledgments:} I would like to thank J. Alev, T. Bridgeland, B. Keller, S. Mozgovoy, R. Pandharipande, B. Siebert, Y. Soibelman and V. Toledano-Laredo for interesting discussions concerning the material in this paper. I am indebted to Y. Soibelman for providing me with a preliminary version of \cite{KS}, and for pointing out a gap in an earlier version of this paper.

\section{Statement of the result}\label{statementresult}

Let $Q$ be a finite quiver, given by a finite set of vertices $I$ and a finite set of arrows, written as $\alpha:i\rightarrow j$ for $i,j\in I$. We assume throughout the paper that $Q$ has no oriented cycles, thus we can order the vertices as $I=\{i_1,\ldots,i_r\}$ in such a way that $k>l$ provided there exists an arrow $i_k\rightarrow i_l$.
Set $\Lambda={\bf Z}I$ with natural basis $i\in I$, and consider the sublattice $\Lambda^+={\bf N}I$.
Let $\Theta\in\Lambda^*$ be a functional (called a stability), and define the slope of a non-zero $d=\sum_{i\in I}d_ii\in\Lambda^+$ as $$\mu(d)=\frac{\Theta(d)}{\dim d},$$
where $\dim d=\sum_{i\in I}d_i$. For $\mu\in{\bf Q}$, define $\Lambda^+_\mu$ as the subsemigroup of $\Lambda^+$ of all $0\not=d\in\Lambda^+$ such that $\mu(d)=\mu$, together with $0\in\Lambda$.\\[1ex]
Define a bilinear form $\langle,\rangle$, the Euler form, on $\Lambda$ via
$$\langle d,e\rangle=\sum_{i\in I}d_ie_i-\sum_{\alpha:i\rightarrow j}d_ie_j\mbox{ for }d,e\in\Lambda.$$
\\[1ex]
Denote by $\{,\}$ the skew-symmetrization of $\langle ,\rangle$, thus $\{d,e\}=\langle d,e\rangle-\langle e,d\rangle$. Define $b_{ij}=\{i,j\}$ for $i,j\in I$.\\[2ex]
We consider the formal power series ring $B={\bf Q}[[\Lambda^+]]={\bf Q}[[x_i: i\in I]]$ with topological basis $x^d=\prod_{i\in I}x_i^{d_i}$ for $d\in\Lambda^+$. The algebra $B$ becomes a Poisson algebra via the Poisson bracket
$$\{x_i,x_j\}=b_{ij}x_ix_j\mbox{ for }i,j\in I.$$
Define automorphisms $T_d$ of $B$ by
$$T_d(x_j)=x_j\cdot(1+x^d)^{\{d,j\}}$$
for all $d\in\Lambda^+$ and $j\in I$. A direct calculation shows that these are indeed Poisson automorphisms of $B$ (see Lemma \ref{pi} for a conceptual proof).\\[1ex]
We study a factorization property in the group ${\rm Aut}(B)$ of Poisson automorphisms of $B$ involving a descending product $\prod^\leftarrow_{\mu\in{\bf Q}}$ indexed by rational numbers. In the following sections, we will see that we actually work in a subgroup of ${\rm Aut}(B)$ where such products are well-defined a priori (see the remark following Definition \ref{prodwell}).

\begin{theorem}\label{mainth} In the group ${\rm Aut}(B)$, we have a factorization
$$T_{i_1}\circ\ldots\circ T_{i_r}=\prod^\leftarrow_{{\mu\in{\bf Q}}}T_\mu,$$
where $$T_\mu(x_j)=x_j\cdot\prod_{i\in I}Q_\mu^i(x)^{b_{ij}}$$ for formal series $Q_\mu^i(x)\in{\bf Z}[[\Lambda^+_\mu]]$. These series $Q_\mu^i(x)$ are given by the generating function
$$Q_\mu^i(x)=\sum_{d\in\Lambda^+_\mu}\chi(M_{d,i}^\Theta(Q))x^d,$$
where $\chi(M_{d,i}^\Theta(Q))$ denotes the Euler characteristic in singular cohomology of a framed moduli space of semistable representations of $Q$ of dimension vector $d$ (see section \ref{smoothmodels} for the precise definition).
\end{theorem}

\section{Quantization}\label{qd}

In this section, we quantize the Poisson algebra $B$ to a skew formal power series ring $A$. We define a class of invertible  elements of $A$ which induce well-defined non-trival automorphisms of $B$ via conjugation.

\begin{definition}\label{sfpsr} For a commutative ring $S$ and an invertible element $q\in S^*$, define a skew formal power series ring $S_q[[\Lambda^+]]$ as follows: as an $S$-module, $S_q[[\Lambda^+]]$ has a topological basis consisting of elements $t^d$ for $d\in\Lambda^+$; multiplication is defined by $$t^d\cdot t^e=q^{-\langle e,d\rangle}t^{d+e}.$$
\end{definition}

We apply this definition to the base ring $K={\bf Q}(q)$, the field of rational functions in $q$, and to $R={\bf Z}[q,q^{-1}]$. This yields $A={\bf Q}(q)_q[[\Lambda^+]]$, with a natural $R$-sublattice $\mathcal{A}=R_q[[\Lambda^+]]$.

\begin{lemma} Specialization at $q=1$ identifies $B$ with $\mathcal{A}/(q-1)\mathcal{A}$ as a Poisson algebra.
\end{lemma}

\proof That specialization at $q=1$ induces an isomorphism of ${\bf Q}$-algebra follows from the definitions. Moreover, we have
$$\{t^d,t^e\}=\frac{t^dt^e-t^et^d}{q-1}=\frac{q^{-\langle e,d\rangle}-q^{-\langle d,e\rangle}}{q-1}t^{d+e}=\{d,e\}t^{d+e}\bmod(q-1).$$\hb

For a series $P(t)=\sum_da_d(q)t^d\in\mathcal{A}$, we denote by $\overline{P}(x)$ its specialization at $q=1$ in $B$, that is, $$\overline{P}(x)=\sum_da_d(1)x^d.$$
Any $R$-algebra automorphism of $\mathcal{A}$ induces a Poisson automorphism of $B$ via specialization; we define a class of automorphisms given by conjugation with invertible elements of $A$ mapping $\mathcal{A}$ to itself.\\[1ex]
For a functional $\eta\in\Lambda^*$, define a twisted form $P(q^\eta t)\in A$ of $P(t)$ by $$P(q^\eta t)=\sum_{d\in\Lambda^+}q^{\eta(d)}a_dt^d.$$
We view any $n\in\Lambda$ as a functional $n\cdot\in\Lambda^*$ by $n\cdot d=\sum_{i\in I}n_id_i$.
A series $P(t)\in A$ with $P(0)=1$, i.e. with constant term equal to $1$, is invertible in $A$.

\begin{proposition}\label{algebra} Let $P(t)\in A$ be a series with constant term equal to $1$, and define $Q^\eta(t)=P(q^\eta t)P(t)^{-1}$ for all $\eta\in\Lambda^*$.
\begin{enumerate}
\item\label{parta} The following conditions on $P(t)$ are equivalent:
\begin{enumerate}
\item $Q^\eta\in\mathcal{A}$ for all $\eta\in\Lambda^*$,
\item $Q^{i\cdot}\in\mathcal{A}$ for all $i\in I$.
\end{enumerate}
\item\label{partb} Conjugation by $P(t)$ maps $t^d$ to $t^d\cdot Q^{\{\_,d\}}(t)$ for all $d\in\Lambda^+$.
\item\label{partc} If the conditions of part \ref{parta} are fulfilled, conjugation by $P(t)$ induces the following Poisson automorphism of $B$:
$$x^d\mapsto x^d\cdot\overline{Q^{\{\_,d\}}}(x).$$
\item\label{partcbis} We have $\overline{Q^\eta}(x)=\prod_{i\in I}\overline{Q^{i\cdot}}(x)^{\eta(i)}$ for all $\eta\in\Lambda^*$.
\item\label{partd} The set of all $P(t)\in A$ with $P(0)=1$ fulfilling the conditions of part \ref{parta} forms a subgroup $\mathcal{S}$ of $A^*$.
\end{enumerate}
\end{proposition}

\proof We have $$Q^{\eta+\nu}(t)=P(q^{\eta+\nu}t)P(t)^{-1}=P(q^{\eta+\nu}t)P(q^\nu t)^{-1}P(q^\nu t)P(t)^{-1}=Q^\eta(q^\nu t)Q^\nu(t).$$
We have $Q^0(t)=1$ by the definitions, so $1=Q^{\eta+(-\eta)}(t)=Q^\eta(q^{-\eta}t)Q^{-\eta}(t)$ and thus $$Q^{-\eta}(t)=Q^\eta(q^{-\eta}t)^{-1}.$$
We conclude that every $Q^\eta(t)$ can be expressed as a product of twisted forms of the $Q^{i\cdot}(t)$ for $i\in I$, proving part \ref{parta} and part \ref{partcbis}.\\[1ex]
We have $t^et^d=q^{\{e,d\}}t^dt^e$ by definition of $A$; it follows that
$$P(t)t^d=\sum_ea_e(q)t^et^d=t^d\sum_eq^{\{e,d\}}a_e(q)t^e=t^dP(q^{\{\_,d\}}t),$$
thus conjugation by $P(t)$ maps $t^d$ to $t^dP(q^{\{\_,d\}}t)P(t)^{-1}$, proving part \ref{partb}. Now part \ref{partc} follows.\\[1ex]
To prove that $\mathcal{S}$ is a subgroup, suppose we are given $P_1(t),P_2(t)\in\mathcal{S}$, with which we associate elements $Q_1^\eta(t),Q_2^\eta(t)\in\mathcal{A}$ as above. We write $Q_2^\eta(t)=\sum_eb_e(q)t^e$. Then
\begin{eqnarray*}
P_1(q^\eta t)P_2(q^\eta t)P_2(t)^{-1}P_1(t)^{-1}&=&P_1(q^\eta t)Q_2^\eta(t)P_1(t)^{-1}\\
&=&P_1(q^\eta t)\sum_eb_e(q)t^e P_1(t)^{-1}\\
&=&\sum_eb_e(q)t^eP_1(q^{\eta+\{\_,e\}}t)P_1(t)^{-1}\\
&=&\sum_eb_e(q)t^eQ_1^{\eta+\{\_,e\}}(t)\end{eqnarray*}
by the identities above. Now every summand in the last sum has coefficients in $R$, thus the sum belongs to $\mathcal{A}$. For $P\in\mathcal{S}$, write $P(q^\eta t)^{-1}P(t)=\sum_dc_d(q)t^d$. We conjugate this equation by $P(t)$ and use part \ref{partb}:
$$Q^\eta(t)^{-1}=P(t)P(q^\eta t)^{-1}=P(t)P(q^\eta t)^{-1}P(t)P(t)^{-1}=\sum_dc_d(q)t^d\underbrace{Q^{\{\_,d\}}(t)}_{\in\mathcal{A}}.$$
The leftmost term belonging to $\mathcal{A}$, we see inductively that all $c_d(q)$ belong to $R$.
This proves part \ref{partd}.\hb

\begin{corollary}\label{cquant} With notation of the previous proposition, the map sending $P(t)\in\mathcal{S}$ to the automorphism
$$x^d\mapsto x^d\cdot\overline{Q^{\{\_,d\}}}(x)$$
defines a group homomorphism $\Phi:\mathcal{S}\rightarrow{\rm Aut}(B)$.
\end{corollary}

As a first example of series in $\mathcal{S}$, we choose a dimension vector $d\in\Lambda^+$ such that $\langle d,d\rangle=1$ (a real root for the root system associated with $Q$), and consider the series
$$P_d(t)=\sum_{n=0}^\infty\frac{q^{-n^2}}{(1-q^{-1})\cdot\ldots\cdot(1-q^{-n})}t^{nd}\in A.$$
The following can be proved using standard identities involving the $q$-binomial coefficients $\left[M\atop N\right]=\frac{(q^M-1)\cdot\ldots\cdot(q^{M-N+1}-1)}{(q^N-1)\cdot\ldots\cdot(q-1)}$. We will give a more conceptual proof in Remark \ref{conceptual}.

\begin{lemma}\label{pi} For any $\eta\in\Lambda$, we have
$$Q_d^\eta(t)=P_d(q^\eta t)P_d(t)^{-1}=\sum_{n=0}^\infty\left[{\eta(d)\atop n}\right]t^{nd}.$$
\end{lemma}

Using Corollary \ref{cquant}, we get in the notation of section \ref{statementresult}:
$$\Phi(P_d(t))=T_d\in{\rm Aut}(B).$$

For later reference, we note the following property:

\begin{lemma}\label{nnn} Suppose $P(t)\in\mathcal{S}$ belongs to ${\bf Q}(q)_q[[\Lambda^+_\mu]]$ for some $\mu\in{\bf Q}$, and let $Q^\eta(t)$ be as in Proposition \ref{algebra}. Then $Q^\eta(t)$ and $Q^\nu(t)$ coincide if $\eta-\nu$ is a rational multiple of the functional $\Theta-\mu\cdot\dim$.
\end{lemma}

\proof The subsemigroup $\Lambda_\mu^+$ has the defining condition $\mu(d)=\mu$, that is, $(\Theta-\mu\cdot\dim)(d)=0$. Thus, the condition in the statement of the lemma is equivalent to $\eta(d)=\nu(d)$ for all $d\in\Lambda_\mu^+$. In the definition of $Q^\eta(t)$, the functional $\eta$ only enters through its values on $\Lambda^+_\mu$ by the choice of $P(t)$; the lemma follows.\hb

\section{Hall algebras and the Harder-Na\-ra\-sim\-han recursion}\label{hallhn}

With a quiver $Q$ and a stability $\Theta$ as before, we associate a system of rational functions defined recursively, and relate it to the cohomology of quiver moduli via Hall algebras; we adapt material of \cite{HNS} to the present setup. Generating series for this system of rational functions will yield the automorphisms $T_\mu$ of Theorem \ref{mainth} via the map $\Phi$ of Corollary \ref{cquant}.

\subsection{The Harder-Narasimhan recursion}

\begin{definition}\label{pdq} Define the following rational functions and their generating series:
\begin{enumerate}
\item for $d\in\Lambda^+$, define $e_d(q)\in{\bf Q}(q)$ by
$$e_d(q)=q^{-\langle d,d\rangle}\prod_{i\in I}\prod_{j=1}^{d_i}(1-q^{-j})^{-1}.$$
\item Define $p_d(q)\in{\bf Q}(q)$ for $d\in\Lambda^+$ recursively as follows: if $\Theta$ is constant on ${\rm supp}(d)=\{i\in I\, :\, d_i\not=0\}$, then $p_d(q)=e_d(q)$. Otherwise, define
$$p_d(q)=e_d(q)-\sum_{d^*}q^{-\sum_{k<l}\langle d^l,d^k\rangle}p_{d^1}(q)\cdot\ldots\cdot p_{d^s}(q),$$
where the sum runs over all non-trivial decompositions $d=d^1+\ldots+d^s$ (i.e. $s\geq 2$ and $d^k\not=0$ for all $k$) such that $\mu(d^1)>\ldots>\mu(d^s)$.
\item Define $P(t)=\sum_{d\in\Lambda^+}e_d(q)t^d\in A$ and $P_\mu(t)=\sum_{d\in\Lambda^+_\mu}p_d(q)t^d\in A$ for all $\mu\in{\bf Q}$.
\end{enumerate}
\end{definition}

\remark The following explicit formula for $p_d(q)$ (a resolution of the defining recursion) is proved in \cite{HNS}:
$$p_d(q)=\sum_{d^*}(-1)^{s-1}q^{-\sum_{k\leq l}\langle d^l,d^k\rangle}\prod_{k=1}^s\prod_{i\in I}\prod_{j=1}^{d_i}(1-q^{-j})^{-1}\in{\bf Q}(q),$$
where the sum ranges over all tuples $d^*=(d^1,\ldots,d^s)$ of non-zero dimension vectors such that
$d=d^1+\ldots+d^s$ and $\mu(d^1+\ldots+d^k)>\mu(d)$ for all $k<s$.

\begin{definition}\label{prodwell} Given elements $c_\mu\in {\bf Q}(q)_q[[\Lambda^+_\mu]]$ for $\mu\in{\bf Q}$ such that $c_\mu(0)=0$, we define
$$\prod^\leftarrow_{\mu\in{\bf Q}}(1+c_\mu)=\sum_{\mu_1>\ldots>\mu_s}c_{\mu_1}\cdot\ldots\cdot c_{\mu_s}.$$
\end{definition}

\remark The sum on the right hand side is indeed well-defined, since calculation of each $t^d$-coefficient reduces to a finite sum. Applying this definition to series $R_\mu\in\mathcal{S}\cap{\bf Q}(q)_q[[\Lambda^+_\mu]]$, we see that decreasing products $\prod^\leftarrow_{\mu\in{\bf Q}}$ in the image of $\Phi:\mathcal{S}\rightarrow{\rm Aut}(B)$ are well-defined via
$$\prod_{\mu\in{\bf Q}}^\leftarrow\Phi(R_\mu)=\Phi(\prod^\leftarrow_{\mu\in{\bf Q}}R_\mu).$$

\begin{lemma}\label{hnsa} We have $P_{i_1}(t)\cdot\ldots\cdot P_{i_r}(t)=P(t)=\prod^\leftarrow_{\mu\in{\bf Q}}P_\mu(t)$ in $A$.
\end{lemma}

\proof We first prove the second identity. By the definition of the functions $p_d(q)$, we have
\begin{eqnarray*}\prod^\leftarrow_{\mu\in{\bf Q}}P_\mu(t)&=&\sum_{\mu_1>\ldots>\mu_s}(\sum_{d\in\Lambda^+_{\mu_1}\setminus 0}p_d(q)t^d)\cdot\ldots\cdot(\sum_{d\in\Lambda^+_{\mu_s}\setminus 0}p_d(q)t^d)\\
&=&\sum_{{{(d^1,\ldots,d^s)}\atop{\mu(d^1)>\ldots>\mu(d^s)}}}p_{d^1}(q)t^{d^1}\cdot\ldots\cdot p_{d^s}(q)t^{d^s}\\
&=&\sum_{{{(d^1,\ldots,d^s)}\atop{\mu(d^1)>\ldots>\mu(d^s)}}}q^{-\sum_{k<l}\langle d^l,d^k\rangle}p_{d^1}(q)\cdot\ldots\cdot p_{d^s}(q)t^{d^1+\ldots+d^s}\\
&=&\sum_d e_d(q)t^d=P_d(t).
\end{eqnarray*}
The first identity follows from the definition of $A$ and the choice of the ordering of the vertices of $Q$; we have
\begin{eqnarray*}P_{i_1}(t)\cdot\ldots\cdot P_{i_r}(t)&=&\sum_{d_1,\ldots,d_r\geq 0}q^{-\sum_{i\in I}d_i^2}\prod_{i\in I}\prod_{j=1}^{d_i}(1-q^{-j})^{-1}t^{d_1i_1}\cdot\ldots\cdot t^{d_ri_r}\\
&=&\sum_{d\in\Lambda^+}q^{\sum_{\alpha:i\rightarrow j}d_id_j-\sum_id_i^2}\prod_{i\in I}\prod_{j=1}^{d_i}(1-q^{-j})^{-1}t^d=P(t).
\end{eqnarray*}\hb

\subsection{Quiver representations}\label{qr}

Let $k$ be a field. We consider finite-dimensional $k$-representations $$M=((M_i)_{i\in I},(M_\alpha:M_i\rightarrow M_j)_{\alpha:i\rightarrow j})$$ of $Q$, given by an $I$-tuple of finite-dimensional $k$-vector spaces $M_i$, together with $k$-linear maps $M_\alpha:M_i\rightarrow M_j$ indexed by the arrows $\alpha:i\rightarrow j$. Let ${\rm mod}_kQ$ be the abelian $k$-linear category of finite dimensional $k$-representations of $Q$. The Grothendieck group of ${\rm mod}_kQ$ can be identified with the lattice $\Lambda$ via the map attaching to the representation $M$ its dimension vector $$\dimv M=\sum_i(\dim_kM_i)i.$$
The above bilinear form $\langle\_,\_\rangle$ then becomes the homological Euler form on ${\rm mod}_kQ$, in the sense that
$$\langle\dimv M,\dimv N\rangle=\dim{\rm Hom}_Q(M,N)-\dim{\rm Ext}^1_Q(M,N)$$
for all representations $M,N$.\\[1ex]
We denote by $S_i$ the one-dimensional simple representation supported at the vertex $i\in I$, and by $P_i$ its projective cover. Then every projective representation of $Q$ is isomorphic to $P^{(n)}=\bigoplus_{i\in I}P_i^{n_i}$ for some $n\in\Lambda^+$.\\[1ex]
For the following basic notions and facts on (semi-)stability of quiver representations, see e.g. \cite{RICRA12}. 
For a non-zero representation $M$, we define its slope as the slope of its dimension vector, i.e. $\mu(M)=\mu(\dimv M)$. We call $M$ semistable if $\mu(U)\leq\mu(M)$ for all non-zero proper subrepresentations $U$ of $M$, and we call $M$ stable if $\mu(U)<\mu(M)$ for all such $U$. Moreover, we call $M$ polystable if it is isomorphic to a direct sum of stable representations of the same slope.\\[1ex]
The semistable representations of a fixed slope $\mu\in{\bf Q}$ form a full abelian subcategory ${\rm mod}^\mu_kQ$, whose simple (resp. semisimple) objects are given by the stable (resp. polystable) representations of slope $\mu$.\\[1ex]
For every representation $M$, there exists a unique Harder-Narasimhan filtration, by which we mean a filtration
$$0=M_0\subset M_1\subset \ldots\subset M_s=M$$
such that all subfactors $M_i/M_{i-1}$ are semistable, and
$$\mu(M_1/M_0)>\ldots>\mu(M_s/M_{s-1}).$$

\subsection{Hall algebras}\label{hallalgebras}

Let $k$ be a finite field. For $d\in\Lambda^+$, fix $k$-vector spaces $M_i$ of dimension $d_i$ for $i\in I$, and let $$R_d=\bigoplus_{\alpha:i\rightarrow j}{\rm Hom}_k(M_i,M_j)$$ be the space of all $k$-representations of $Q$ on the vector spaces $M_i$, on which the group $$G_d=\prod_{i\in I}{\rm GL}(M_i)$$ acts via base change $$(g_i)_i\cdot(M_\alpha)_\alpha=(g_jM_\alpha g_i^{-1})_{\alpha:i\rightarrow j},$$
such that the $G_d$-orbits in $R_d$ correspond naturally to the isomorphism classes of representations of $Q$. Let ${\bf Q}^{G_d}(R_d)$ be the space of (arbitrary) $G_d$-invariant ${\bf Q}$-valued functions on $R_d$, and define $$H_k((Q))=\prod_{d\in\Lambda^+}{\bf Q}^{G_d}(R_d).$$
This becomes a $\Lambda^+$-graded algebra, the (completed) Hall algebra of $Q$, via the convolution type product
$$(fg)(M)=\sum_{U\subset M}f(U)g(M/U)$$
for functions $f\in{\bf Q}^{G_d}(R_d)$, $g\in{\bf Q}^{G_e}(R_e)$ and representations $M\in R_{d+e}$. Note that the sum over all subrepresentations $U$ is finite and that the value of $f$ on $U$ (resp. of $g$ on $M/U$) is well-defined by the definitions.\\[1ex]
Let ${\bf Q}_{|k|}[[\Lambda^+]]$ be defined as in Definition \ref{sfpsr}, where the cardinality $|k|$ of $k$ is viewed as an element of ${\bf Q}$. The map
$$\int:H((Q))\rightarrow{\bf Q}_{|k|}[[\Lambda^+]]$$
given by
$$\int f=\frac{1}{|{G_d}|}\sum_{M\in R_d}f(M)t^d$$
for $f\in{\bf Q}^{G_d}(R_d)$ is a ${\bf Q}$-algebra morphism by \cite{HNS}.\\[1ex]
Define $1_d=1_{R_d}$ as the characteristic function of $R_d$, and define $1_d^{sst}=1_{R_d^{sst}}$ as the characteristic function of the locus $R_d^{sst}$ of semistable representations in $R_d$.
We form generating functions of these elements by $$e:=\sum_{d\in\Lambda^+}1_d,\;\;\; e_\mu=\sum_{d\in\Lambda^+_\mu}1_d^{sst}
\mbox{ for }\mu\in{\bf Q}.$$

\begin{lemma}\label{hns1} For every $d\in\Lambda^+$, we have the following identity in $H_k((Q))$:
$$1_d=\sum_{d^*}1_{d^1}^{sst}\cdot\ldots\cdot 1_{d^s}^{sst},$$
the sum running over all decompositions $d=d^1+\ldots+d^s$ into non-zero dimension vectors such that $\mu(d^1)>\ldots>\mu(d^s)$. Consequently, we have $e=\prod^\leftarrow_{\mu\in{\bf Q}}e_\mu$ in $H_k((Q))$.
\end{lemma}

\proof The existence and uniqueness of the Harder-Narasimhan filtration (see section \ref{qr}) can be rephrased as follows using the definition of the Hall algebra: for every $k$-representation $M$ of $Q$, there exists a unique tuple $(d^1,\ldots,d^s)$ of dimension vectors such that $\mu(d^1)>\ldots>\mu(d^s)$ and 
$$(1_{d^1}^{sst}\cdot\ldots\cdot 1_{d^s}^{sst})(M)=1$$
(namely, the $d^i$ are the dimension vectors of the subfactors in the Harder-Narasimhan filtration). The first identity follows. Similar to the proof of Lemma \ref{hnsa} above, this allows us to compute
\begin{eqnarray*}\prod^\leftarrow_{\mu\in{\bf Q}}e_\mu&=&\sum_{\mu_1>\ldots>\mu_s}(\sum_{d\in\Lambda^+_{\mu_1}\setminus 0}1_d^{sst})\cdot\ldots\cdot(\sum_{d\in\Lambda^+_{\mu_s}\setminus 0}1_d^{sst})\\
&=&\sum_{{{(d^1,\ldots,d^s)}\atop{\mu(d^1)>\ldots>\mu(d^s)}}}1_{d^1}^{sst}\cdot\ldots\cdot 1_{d^s}^{sst}\\
&=&\sum_{d\in\Lambda^+}1_d=e_d.\end{eqnarray*}\hb

Since the rational functions $e_d(q)$ and $p_d(q)$ have no poles at $q=|k|$ by Definition \ref{pdq}, we can specialize the generating functions $P(t)$ and $P_\mu(t)$ to ${\bf Q}_{|k|}[[\Lambda^+]]$.

\begin{proposition}\label{spec} The series $P(t)$ specializes to $\int e$, and the series $P_\mu(t)$ specialize to $\int e_\mu$ for all $\mu\in{\bf Q}$.
\end{proposition}

\proof By the definitions of $R_d$, $G_d$ and the Euler form $\langle\_,\_\rangle$ on $Q$, we have
$$\int 1_d=\frac{|R_d|}{|G_d|}=
\frac{
{|k|}^{\sum_{\alpha:i\rightarrow j}d_id_j}
}
{
\prod_{i\in I}|{\rm GL}_{d_i}(k)|
}=|k|^{-\langle d,d\rangle}\prod_{i\in I}\prod_{j=1}^{d_i}(1-|k|^{-j})^{-1}=e_d(|k|),$$
proving the first statement. The second now follows since the elements $\int 1_d^{sst}$ satisfy the same recursion as the elements $p_d(q)$ by Lemma \ref{hns1}.\hb

\remark By the definition of the function $1_d^{sst}$, we thus have 
$$p_d(|k|)=\frac{|R_d^{sst}|}{|G_d|}.$$
In case $d$ is coprime for $\Theta$, this is used in \cite{HNS} to prove that
$$(q-1)\cdot p_d(q)=\sum_i\dim H^i(M_d^{st}(Q),{\bf Q})q^{i/2},$$
where $M_d^{st}(Q)$ denotes the moduli space of stable representations of $Q$ of dimension vector $d$ (see section \ref{smoothmodels} for the definitions). The result holds since in the coprime case, we have a smooth projective moduli space whose numbers of rational points over finite fields $k$ behave polynomially in $|k|$.

\section{Smooth models and the proof of Theorem \ref{mainth}}\label{smpt}

Using the construction of (framed versions of) moduli spaces of representations of quivers, we prove that the series $P_\mu$ belong to $\mathcal{S}$. This fact, together with the Harder-Narasimhan recursion of the previous section, proves Theorem \ref{mainth}.

\subsection{Quiver moduli}\label{smoothmodels}

In this section, we work over the complex numbers. For every dimension vector $d\in\Lambda^+$, there exists a smooth complex variety $M_d^{st}(Q)$ whose points parametrize isomorphim classes of stable complex representations of $Q$ of dimension vector $d$. It embeds as an open subset into a (typically singular) projective variety $M_d^{sst}(Q)$ whose points parametrize isomorphism classes of polystable representations of $Q$ of dimension vector $d$. If $d$ is coprime for $\Theta$, by which we mean that $\mu(e)\not=\mu(d)$ for all $0\leq e<d$, we have $M_d^{st}(Q)=M_d^{sst}(Q)$, consequently a smooth projective complex variety.\\[1ex]
For $n\in\Lambda^+$, fix additional vector spaces $V_i$ of dimension $n_i$ for $i\in I$.

\begin{theorem}[\cite{SM}]\label{defsm} There exists a projective variety $M_{d,n}(Q)$ parametrizing equivalence classes of pairs $(M,f)$ consisting of a semistable representation $M$ of $Q$ on the vector spaces $M_i$, together with a tuple of maps $f=(f_i:V_i\rightarrow M_i)_{i\in I}$
such that $\mu(U)<\mu(M)$ whenever $U\subset M$ is a proper subrepresentation of $M$ containing the image of $f$ (that is, the subrepresentation generated by the $f_i(V_i)$);
such pairs are considered up to isomorphisms of the representations intertwining the additional maps, that is, $(M,f)$ is equivalent to $(M',f')$ if there exists an isomorphism $\varphi:M\rightarrow M'$ such that $f'=\varphi f$.
\end{theorem}

This is a framed version of moduli of (semistable) quiver representations, called smooth models in \cite{SM}. There it is also shown that the framing datum $f$ induces a morphism from the projective representation $P^{(n)}$ to $M$ such that $\mu(U)<\mu(M)$ for any proper subrepresentation $U$ of $M$ containing its image.\\[1ex]
There exists a canonical projective morphism $\pi:M_{d,n}^\Theta(Q)\rightarrow M_d^{sst}(Q)$, which by \cite{SM} is \'etale locally trivial for a suitable (Luna type) stratification of $M_d^{sst}(Q)$ with known fibres (they are isomorphic to certain nilpotent parts of smooth models for the trivial stability and quivers with oriented cycles).\\[1ex]
As a special case, we consider the quiver $Q^0$ consisting of a single vertex and no arrows, with trivial stability. From the definitions, it is easy to see that $M_{d,n}^0(Q^0)\simeq{\rm Gr}_d^n({\bf C})$, the Grassmannian of $k$-planes in $n$-space.

\subsection{Proof of Theorem \ref{mainth}}

Let again $k$ be a finite field. The semistable representations of $Q$ of fixed slope $\mu\in{\bf Q}$ form a full abelian subcategory ${\rm mod}_k^\mu(Q)$ of ${\rm mod}_k^\mu(Q)$. Thus, for a subrepresentation $U\subset M$ of a representation $M\in{\rm mod}_k^\mu(Q)$, we can form the intersection
$$\langle U\rangle_\mu=\bigcap_{{V\in{\rm mod}_k^\mu(Q)}\atop{U\subset V\subset M}}V\in{\rm mod}_k^\mu(Q).$$
This is the minimal subrepresentation of $M$ containing $U$ which is semistable of slope $\mu$. For $n\in\Lambda^+$, we denote by ${\rm Hom}^0_Q(P^{(n)},M)$ the set of all morphisms $f\in{\rm Hom}_Q(P^{(n)},M)$ such that $\langle{\rm Im}(f)\rangle_\mu=M$. As mentioned above, the data $(M,f)$ of Theorem \ref{defsm} can be viewed as pairs consisting of a representation $M\in{\rm mod}_k^\mu(Q)$, together with a map $f\in{\rm Hom}_Q^0(P^{(n)},M)$.\\[1ex]
For $d\in\Lambda^+_\mu$, define $1_{d,n}^{sst}$ as the function taking value $|{\rm Hom}_Q(P^{(n)},M)|=q^{n\cdot d}$ on $M\in R_d^{sst}$, and value $0$ outside $R_d^{sst}$. Define $f_{d,n}\in H_k((Q))$ as the function taking value $|{\rm Hom}_Q^0(P^{(n)},M)|$ on $M\in R_d^{sst}$, and value $0$ outside $R_d^{sst}$. We form generating functions $e_{\mu,n}=\sum_{d\in\Lambda^+_\mu}1_{d,n}^{sst}$ and $h_{\mu,n}=\sum_{d\in\Lambda^+_\mu}f_{d,n}$.

The following lemma is a special case of \cite[Lemma 5.1]{SM}:

\begin{lemma} We have $e_{\mu,n}=h_{\mu,n}\cdot e_\mu$ in $H_k((Q))$.
\end{lemma}

\proof For a representation $M\in{\rm mod}_k^\mu(Q)$, we have by the definitions
$$(h_{\mu,n}\cdot e_\mu)(M)=\sum_U|{\rm Hom}^0(P^{(n)},U)|,$$
where the sum runs over all subrepresentations $U$ of $M$ such that $U\in{\rm mod}_k^\mu(Q)$. But the set of all pairs $(U,f)$ consisting of a subrepresentation $U\subset M$ such that $U\in{\rm mod}_k^\mu(Q)$ and a morphism $f\in{\rm Hom}^0(P^{(n)},U)$ is naturally in bijection to ${\rm Hom}(P^{(n)},M)$, by associating to $f:P^{(n)}\rightarrow M$ the pair $(\langle {\rm Im}(f)\rangle_\mu,f:P^{(n)}\rightarrow\langle{\rm Im}(f)\rangle_\mu)$. Thus,
$$(h_{\mu,n}\cdot e_\mu)(M)=|{\rm Hom}(P^{(n)},M)|=e_{\mu,n}(M),$$
proving the lemma.\hb

\begin{proposition}\label{proppmu} The series $P_\mu(t)$ belongs to $\mathcal{S}$. The automorphism $\Phi(P_\mu(t))$ is given by
$$x^d\mapsto x^d\cdot\prod_{i\in I}Q_\mu^i(x)^{\{i,d\}},$$
where
$$Q_\mu^i(x)=\sum_{d\in\Lambda_d^+}\chi(M_{d,i}^\Theta(Q))x^d.$$
\end{proposition}

\proof By Proposition \ref{spec}, we know that the series $P_\mu(t)$ specializes to $\int e_\mu$ under specialization of $q$ to $k$. By definition of $e_{\mu,n}$, the series $P_\mu(q^{n\cdot}t)$ thus specializes to $\int e_{\mu,n}$ for $n\in\Lambda^+$. By the previous lemma, the function $Q_\mu^{n\cdot}(t)=P_\mu(q^{n\cdot}t)P_\mu(t)^{-1}$ thus specializes to $\int e_{\mu,n}^0$. By definition, the function $f_{d,n}$ integrates to 
$$\int f_{d,n}=\frac{|\{(M,f):M\in R_d^{sst},\; f\in{\rm Hom}^0(P^{(n)},M\}|}{|G_d|}.$$
As in the proof of \cite[Theorem 5.2]{SM}, we can conclude, using the definition of the smooth model $M_{d,n}^\Theta(Q)$ and the comparison between numbers of rational points and Betti numbers as in the final remark of section \ref{hallalgebras}, that the Poincare polynomial $\sum_i\dim H^i(M_{d,n}^\Theta(Q),{\bf Q})q^{i/2}$ specializes to $\int f_{d,n}$ at $q=|k|$. Consequently, we have
$$\int e_{\mu,n}=\sum_{d\in\Lambda^+_\mu}\sum_i\dim H^i(M_{d,n}^\Theta(Q),{\bf Q})q^{i/2}t^d.$$
Thus, the coefficients $b_d(q)$ of $Q_\mu^{n\cdot}(t)$ are rational functions taking integer values at all prime powers $q$, and thus are polynomials in $q$. By part \ref{parta} of Proposition \ref{algebra}, we conclude that $P_\mu(t)\in\mathcal{S}$. Using parts \ref{partc} and \ref{partcbis} of the same proposition, we derive the claimed formula for $\Phi(P_\mu(t))$ since the Poincare polynomial specializes to the Euler characteristic at $q=1$ in absence of odd cohomology. \hb

\remark\label{conceptual} Applying this result to the quiver $Q^0$, we get Lemma \ref{pi} in the case $d=i$. The assumption $\langle d,d\rangle=1$ there allows to generalize to such $d$ using the definition of $P_d(q)$.\\[2ex]
Theorem \ref{mainth} now follows immediately from Lemma \ref{pi}, Lemma \ref{hnsa} and Proposition \ref{proppmu}.

\section{Examples}\label{examples}

\subsection{Generalized Kronecker quivers}

Let $Q=K_m$ be the $m$-arrow Kronecker quiver with set of vertices $I=\{i,j\}$ and $m$ arrows from $j$ to $i$. We have $\{ai+bj,ci+dj\}=m\cdot(ad-bc)$ and in particular $b_{ij}=m$. We choose the stability $\Theta=j^*$. Writing a dimension vector $d\in\Lambda^+$ as $d=ai+bj$ for $a,b\in{\bf N}$, the automorphism $T_{a,b}=T_d$ of the Poisson algebra $B_m={\bf Q}[[x_i,x_j]]$ with Poisson bracket $\{x_i,x_j\}=mx_ix_j$ is then given by 
$$T_{a,b}:\left\{\begin{array}{lll}x_i&\mapsto&x_i(1+x_i^{a}x_j^{b})^{-mb}\\ x_j&\mapsto&x_j(1+x_i^{a}x_j^{b})^{ma}\end{array}\right\}.$$

We prove the following weak version of \cite[Conjecture 1]{KS}:

\begin{theorem} In ${\rm Aut}(B_m)$, there exists a factorization
$$T_i\circ T_j=\prod^\leftarrow_{b/a\in{\bf Q}}T_{a,b}^{d(a,b)},$$
where $d(a,b)\in\frac{1}{{\rm gcd}(a,b)}{\bf Z}$ for all $a,b\in{\bf N}$.
\end{theorem}

\proof By Theorem \ref{mainth}, we have a factorization $T_i\circ T_j=\prod^\leftarrow_{\mu\in{\bf Q}}T_\mu$, where $T_\mu$ is given by
$$T_\mu:\left\{\begin{array}{lll}x_i&\mapsto&x_i\cdot Q_\mu^j(x)^{-m}\\ x_j&\mapsto&x_j\cdot Q_\mu^i(x)^m\end{array}\right\}$$
for $Q_\mu^i(x)$, $Q_\mu^j(x)$ being the corresponding generating series for Euler characteristics of smooth models.
Writing $\mu=b/(a+b)$ for coprime $a,b\in{\bf N}$, we have $\Lambda^+_\mu={\bf N}\cdot(a,b)$, and thus $Q_\mu^i(x),Q_\mu^j(x)\in{\bf Z}[[x_i^ax_j^b]]$.
We choose integers $c,d\in{\bf Z}$ such that $ac+bd=1$ and define $$F_\mu(x)=Q_\mu^i(x)^cQ_\mu^j(x)^d\in{\bf Z}[[x_i^ax_j^b]].$$
By Lemma \ref{nnn}, we have $Q_\mu^j(x)^a=Q_\mu^i(x)^b$. It then follows that $$F_\mu(x)^a=Q_\mu^i(x)\mbox{ and }F_\mu(x)^b=Q_\mu^j(x).$$
We can factor $F_\mu(x)$ into an infinite product $$F_\mu(x)=\prod_{k\geq 1}(1+(x_i^ax_j^b)^k)^{c(\mu,k)}$$
for integers $c(\mu,k)\in{\bf Z}$. These two identities allow us to write $T_\mu$ in the form
$$T_\mu:\left\{\begin{array}{lll}x_i&\mapsto&x_i\cdot\prod_{k\geq 1}(1+x_i^{ka}x_j^{kb})^{-mbc(\mu,k)},\\
x_j&\mapsto&x_j\cdot\prod_{k\geq 1}(1+x_i^{ka}x_j^{kb})^{mac(\mu,k)}\end{array}\right\}.$$
Using the definition of the automorphisms $T_{ka,kb}$ and defining $d(ka,kb)=\frac{1}{k}c(\mu,k)$, the theorem follows.\hb

It is now easy to re-derive the examples of \cite[Section 1.4]{KS}. For $m=1$, there are three isomorphism classes of stable representations of $K_1$, namely the two simple representations $S_i$, $S_j$ and their non-trivial extension $X$. Thus, the relevant smooth models reduce to Grassmannians as in the case of the quiver $Q^0$, yielding the factorization
$$T_{1,0}\circ T_{0,1}=T_{0,1}\circ T_{1,1}\circ T_{0,1},$$
which can also be verified directly.\\[1ex]
For $m=2$, we have a unique stable representation for any dimension vector $ki+(k+1)j$ or $(k+1)i+kj$. Moreover, we have a ${\bf P}^1$-family of stable representations of dimension vector $i+j$, and no stables of dimension vector $ki+kj$ for $k\geq 2$. We have $M_{ki+kj,n}^\Theta(K_2)\simeq{\rm Hilb}^k({\bf P}^1)\simeq{\bf P}^k$ for $n=i,j$, and thus
$$Q_{1/2}^i(x)=Q_{1/2}^j(x)=\sum_{k\geq 0}(k+1)(x_ix_j)^k=(1-x_ix_j)^{-2}.$$
We arrive at the following factorization in $B_2$:
$$T_{1,0}\circ T_{0,1}=(T_{0,1}\circ T_{1,2}\circ T_{2,3}\circ\ldots)\circ T\circ (\ldots\circ T_{3,2}\circ T_{2,1}\circ T_{1,0}),$$
where
$$T:\left\{\begin{array}{lll}x_i&\mapsto&x_i\cdot(1-x_ix_j)^4\\ x_j&\mapsto&x_j\cdot(1-x_ix_j)^{-4}\end{array}\right\}.$$

\remark In \cite[Theorem 0.1]{GPS}, the automorphisms $T_\mu$, given by generating functions of Euler characteristics of quiver moduli, are calculated in terms of generating function of genus $0$ Gromov-Witten invariants of toric surfaces. This alternative interpretation hints towards a link between the underlying geometries.

\subsection{Dynkin quivers} 

Let $Q$ be a quiver of Dynkin type. Then the isomorphism classes of indecomposable representations $U_\alpha$ correspond bijectively to the positive roots $\alpha\in\Delta^+$ of the corresponding root system. The Harder-Narasimhan filtration of section \ref{qr} can be replaced by a more explicit filtration (we refer to \cite{RFeigin} for details and a related application of this filtration). Namely, the positive roots can be ordered as $\alpha_1,\ldots,\alpha_\nu$ in such a way that ${\rm Hom}_Q(U_{\alpha_l},U_{\alpha_k})=0={\rm Ext}^1_Q(U_{\alpha_k},U_{\alpha_l})$ for $k<l$. Now every representation $M$ of $Q$ is of the form $M=\bigoplus_kU_{\alpha_k}^{m_k}$, and thus admits a unique filtration
$$0=M^\nu\subset M^{\nu-1}\subset\ldots\subset M^1\subset M^0=M$$
such that $M^{k-1}/M^{k}\simeq U_{\alpha_k}^{m_k}$ for all $k=1,\ldots,\nu$. The arguments of section \ref{hallhn} can be applied to this filtration, yielding a factorization $P_{i_1}\cdot\ldots\cdot P_{i_r}=P_{\alpha_\nu}\cdot\ldots\cdot P_{\alpha_1}$.
Applying the map $\Phi:\mathcal{S}\rightarrow{\rm Aut}(B)$, we conclude
$$T_{i_1}\circ\ldots\circ T_{i_r}=T_{\alpha_\nu}\circ\ldots\circ T_{\alpha_1}.$$

\end{document}